\documentclass[11pt]{article}
\usepackage{dsfont}
\usepackage{graphics}
\usepackage[demo]{graphicx}
\usepackage{epsfig}
\usepackage{pstricks}
\usepackage[normal]{subfigure}
\usepackage[latin1]{inputenc}
\usepackage[english,francais]{babel}
\usepackage{relsize,exscale}
\usepackage{makeidx}
\usepackage{enumitem}
\usepackage{amsfonts,amssymb,amsmath}
\usepackage{graphicx}
\usepackage{color}
\usepackage{multirow}
\usepackage{mathrsfs}
\usepackage[normalem]{ulem}
\usepackage{cancel}
\usepackage{bbm}
\newenvironment{prooff}{{\it Proof :}}{\hfill\rule{2mm}{2mm}\vskip3mm\par}
\newtheorem{theorem}{Theorem}[section]
\newtheorem{lemma}[theorem]{Lemma}

\newtheorem{e-definition}[theorem]{Definition\rm}
\newtheorem{remark}{\it Remark\/}

\usepackage{color}
\definecolor{dred}{rgb}{0.92,0,0}
\definecolor{dgreen}{rgb}{0,0.92,0}
\definecolor{dblue}{rgb}{0,0,0.92}
\definecolor{dyellow}{rgb}{0.95,0.95,0}

\newcommand{\R}{\mathbb{R}}
\newcommand{\N}{\mathbb{N}}
\def\D{\displaystyle}
\newcommand{\hs}{\hspace{0.1cm}}

\newcommand{\sa}{\\ [0.2cm]}
\usepackage{multirow}
\graphicspath{
{./Figures/}
{./}
}
\title{
A new mixed functional-probabilistic approach\\ for finite element accuracy}
\author{Jo\"el Chaskalovic \thanks{D'Alembert,
Sorbonne University, Paris, France, (\emph{email}: jch1826@gmail.com)}
\qquad
Franck Assous
\thanks{Department of Mathematics,
Ariel University, 40700 Ariel, Isra\"el, (\emph{email}: franckassous55@gmail.com).}
}
\date{}
\begin{document}
\maketitle
\selectlanguage{english}

\begin{abstract}
\noindent The aim of this paper is to provide a new perspective on finite element accuracy. Starting from a geometrical reading of the Bramble-Hilbert lemma, we recall the two probabilistic laws we got in previous works that estimate the relative accuracy, considered as a random variable, between two finite elements $P_k$ and $P_m$, ($k < m$). Then, we analyze the asymptotic relation between these two probabilistic laws when the difference $m-k$ goes to infinity. New insights which qualified the relative accuracy in the case of high order finite elements are correspondingly obtained. \\[0.2cm]
\noindent{\footnotesize {\em keywords}: Error estimates, Finite elements, Bramble-Hilbert lemma, Probability.}
\end{abstract}

\section{Introduction}\label{intro}
\noindent Error estimates play a crucial role in advancement of finite element methods. The development and the actual use of a given numerical method is often, if not always, guided by its performance in terms of reliability and accuracy. For this reason, it is still an active subject of research for a wide range of applied mathematicians.\sa
Since the seminal papers of Strang\&Fix \cite{StFi73}, Ciarlet\&Raviart \cite{Ciarlet_Raviart},  Babuska\cite{Babu71}, Bramble\&Hilbert \cite{BrHi70}, with co-workers, a plethora of work has been published to elaborate and improve the error estimates in various configurations. The goal is to find bounds for the error $u-u_h$, between the exact solution $u$ of a partial differential equation and its finite element approximation $u_h$.\sa
Usually, the main property of error estimates which is considered concerns the rate of convergence of a given finite element. In general, these estimates tell us that the finite element error $\|u-u_h\|$, for a given chosen norm, is $O(h^k)$, where $h$ is the mesh size, namely the largest diameter of the elements in a given mesh, and $k$ a positive integer. As the constant involved in $O(h^k)$ is in most cases unknown, it very seldom considered in the analysis.\sa
The aim of these estimates is generally to give a measure of the efficiency of the considered finite element method, and tell us how fast the error decreases as we decrease the mesh size $h$. However, in these estimates, this constant depends, among others, on the basis functions of the concerned finite element method, and on a semi-norm of the exact solution $u$, (see for instance \cite{Arcangeli_Gout}).\sa
Furthermore, quantitative uncertainties are commonly produced in the mesh generation so that quantitative uncertainties also exist in the approximate solution $u_h$. For this reason, we have considered the approximation error as a random variable \cite{CMAM_2019}, and we aimed to evaluate the probability of the difference between two approximation errors $u- u_h^{(k)}$ and $u- u_h^{(m)}$, for a suited functional norm, corresponding to Lagrange finite elements $P_k$ and $P_m, (k<m)$. It is the reason why we introduced in \cite{CMAM_2019} a probabilistic framework to compare the relative accuracy between these two finite elements.\sa
This paper is mainly devoted to the asymptotic relation between the two probabilistic laws we derived in \cite{CMAM_2019}. Amongst other, it will highlight the relative accuracy between high order finite elements. To the best of our knowledge, this is the first time that such mixed functional and probabilistic approaches are combined to provide new perspectives on finite element accuracy. \sa
The paper is structured as follows: Section \ref{Geo_and_Proba} summarizes the results of \cite{CMAM_2019} necessary for understanding the rest of our analysis: the geometrical interpretation of error estimates and the probabilistic laws we got for the relative finite element accuracy.  In Section \ref{propPk} properties of $P_k$ basis polynomials are derived whereas Section \ref{Asymptotic_limit} is devoted to the asymptotic relationship between the two probability laws of Section \ref{Geo_and_Proba}. Concluding remarks follow.

\section{Geometrical interpretation of error estimates and related probabilistic laws}\label{Geo_and_Proba}
\begin{figure}[h]
  \centering
  \includegraphics[width=10cm]{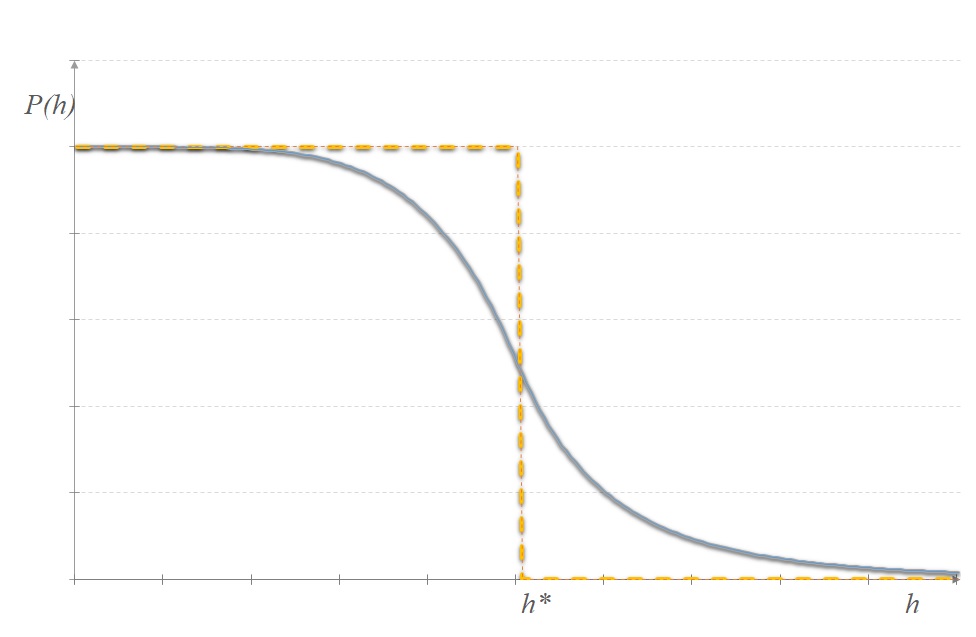}
  \caption{Case $m-k\neq 1$: shape of the sigmoid distribution (\ref{Nonlinear_Prob}) (full line) and the two steps corresponding one (\ref{Heaviside_Prob}) (dashed line), $(P(h)\equiv Prob\{X^{(m)}(h) \leq X^{(k)}(h)\})$.} \label{Sigmoid}
\end{figure}

\noindent We consider an open bounded and non empty subset $\Omega$ of $\R^{n}$, and $\Gamma$ its boundary assumed to be $C^1- $piecewise. Let $u$ be the solution to the second order elliptic variational formulation:
\begin{equation}\label{VP}
\textbf{(VP}\textbf{)} \hspace{0.2cm} \left\{
\begin{array}{l}
\mbox{Find } u \in   V \mbox{ solution to:} \\[0.1cm]
a(u,v) = l(v), \quad\forall v \in V,
\end{array}
\right.
\end{equation}
where $V$ is a given Hilbert space endowed with a norm $ \left\|.\right\|_{V}$, $a(\cdot,\cdot)$ is a bilinear, continuous and $V-$elliptic form defined on $V \times V$, and $l(\cdot)$ a linear continuous form defined on~$V$.\sa
Classically, variational problem \textbf{(VP)} has, one and only one, solution $u \in V$ (see for example \cite{ChaskaPDE}). In this paper, we will restrict ourselves to the case where $V$ is the usual Sobolev space of distributions $H^1(\Omega)$. \sa
Let us also consider an approximation $u_{h}$ of $u$, solution  to the approximate variational formulation:
\begin{equation}\label{VP_h}
\textbf{(VP}\textbf{)}_{h} \hspace{0.2cm} \left\{
\begin{array}{l}
\mbox{Find } u_{h} \in   V_h \mbox{ solution to:} \\[0.1cm]
a(u_{h},v_{h}) = l(v_{h}),\quad \forall v_{h} \in V_h, 
\end{array}
\right.
\end{equation}
where $V_h$ is a given finite-dimensional subset of $V$. \sa
To state a corollary of Bramble-Hilbert's lemma and a corresponding error estimate, we follow \cite{RaTho82} or \cite{Ciarlet}, and we assume that $\Omega$ is exactly covered by a mesh ${\mathcal T}_h$ composed by $N_K$ n-simplexes $K_{\mu}, (1 \leq \mu \leq N_K),$ which respects classical rules of regular discretization, (see for example \cite{ChaskaPDE} for the bidimensional case and \cite{RaTho82} in $\R^n$). Moreover, we denote by $P_k(K_{\mu})$ the space of polynomial functions defined on a given n-simplex $K_{\mu}$ of degree less than or equal to $k$, ($k \geq$ 1). \sa
Then, we remind below the result of \cite{RaTho82} from which our study is developed: \vspace{0.1 cm}
\begin{lemma}\label{Thm_error_estimate}
Suppose that there exists an integer $k \geq 1$ such that the approximation $u_h$ of $V_h$ is a continuous piecewise function composed by polynomials which belong to $P_k(K_{\mu}), (1\leq \mu\leq  N_K)$. \sa
Then, $u_h$ converges to $u$ in $H^1(\Omega)$:
\begin{equation}
\D\lim_{h\rightarrow 0}\|u_h-u\|_{1,\Omega}=0.
\end{equation}
Moreover, if the exact solution $u$ belongs to $H^{k+1}(\Omega)$, we have the following error estimate:
\begin{equation}\label{estimation_error}
\|u_h-u\|_{1,\Omega} \hs \leq \hs \mathscr{C}_k\,h^k \, |u|_{k+1,\Omega}\,,
\end{equation}
where $\mathscr{C}_k$ is a positive constant independent of $h$, $\|.\|_{1,\Omega}$ the classical norm in $H^1(\Omega)$ and $|.|_{k+1,\Omega}$ denotes the semi-norm in $H^{k+1}(\Omega)$.
\end{lemma}
%
Consider now two families of Lagrange finite elements $P_k$ and $P_m$ corresponding to two values $(k,m)\in \N^{*2}$, $(k < m)$, the corresponding inequalities given by (\ref{estimation_error}), assuming that the solution $u$ to \textbf{(VP)} belongs to $H^{m+1}(\Omega)$, are:
\vspace{-0.2cm}
\begin{eqnarray}
\|u^{(k)}_h-u\|_{1,\Omega} \hs & \leq & \hs \mathscr{C}_k h^{k}\, |u|_{k+1,\Omega}, \label{Constante_01} \\
\|u^{(m)}_h\hspace{-0.09cm}-u\|_{1,\Omega} \hs & \leq & \hs \mathscr{C}_m h^{m}\, |u|_{m+1,\Omega}\,, \label{Constante_02}
\end{eqnarray}
where $u^{(k)}_h$ and $u^{(m)}_h$ respectively denotes the $P_k$ and $P_m$ Lagrange finite element approximations of $u$.\\[0.2cm]
Now, if one considers a given mesh for the finite element $P_m$ that contains the mesh processed for the $P_k$ approximation, then for the particular class of problems where $\textbf{(VP)}$ is equivalent to a minimization formulation $\textbf{(MP)}$, (see for example \cite{ChaskaPDE}), one can show that the approximation error for $P_m$ is always lower than the one for $P_k$, and $P_m$ is more accurate than $P_k$, for all values of the mesh size $h$.\sa
In this paper we consider a more general where, for a given mesh size $h$, two independent meshes for $P_k$ and $P_m$ are built by a mesh generator. So, usually, to compare the relative accuracy between these two finite elements, one asymptotically considers inequalities (\ref{Constante_01}) and (\ref{Constante_02}) to conclude that, when $h$ goes to zero, $P_m$ is more accurate that $P_k$, as $h^m$ goes faster to zero than $h^k$. \sa
However, for each application, $h$ has a fixed value and this way of comparison is no longer valid. For this reason, our viewpoint will be to determine the relative accuracy between $P_k$ and $P_m, (k<m)$, for a value of $h$ corresponding to two independent meshes.\sa
To this end, let us set:
\begin{equation}\label{Ck_Cm}
C_k = \mathscr{C}_k |u|_{k+1,\Omega} \mbox{ and } C_m = \mathscr{C}_m |u|_{m+1,\Omega}.
\end{equation}
Therefore, instead of (\ref{Constante_01}) and (\ref{Constante_02}), we consider in the sequel the two following inequalities:
\begin{eqnarray}
\|u^{(k)}_h-u\|_{1,\Omega} \hs & \leq & \hs C_k h^{k}, \label{Constante_01_2} \\
\|u^{(m)}_h\hspace{-0.09cm}-u\|_{1,\Omega} \hs & \leq & \hs C_m h^{m}. \label{Constante_02_2}
\end{eqnarray}
So, we proposed in \cite{CMAM_2019} a geometrical interpretation of (\ref{Constante_01_2})-(\ref{Constante_02_2}) which enabled us to consider the values of $\|u^{(k)}_h-u\|_{1,\Omega}$ and $\|u^{(m)}_h-u\|_{1,\Omega}$ as two random variables, respectively denotes $X^{(k)}$ and $X^{(m)}$, whose values belong to $[0, C_i h^i], (i=k \mbox{ or } i=m),$ due to (\ref{Constante_01_2})-(\ref{Constante_02_2}).\sa
Then, we derived two probabilistic laws of the event:
\begin{equation}
\D\left\{X^{(m)}(h) \leq X^{(k)}(h)\right\} \equiv \left\{\|u^{(m)}_h-u\|_{1,\Omega} \leq \|u^{(k)}_h-u\|_{1,\Omega}\right\},
\end{equation}
which corresponds to the relative accuracy between the two finite elements $P_k$ and $P_m$, for a given value of the mesh size $h$.\sa
More precisely, let us introduce the two random events $A$ and $B$ as follows:
\vspace{-0.2cm}
\begin{eqnarray}
A & \equiv & \left\{\|u^{(m)}_h-u\|_{1,\Omega} \leq \|u^{(k)}_h-u\|_{1,\Omega} \right\}, \label{A}\\
B & \equiv & \left\{\|u^{(k)}_h-u\|_{1,\Omega}\in [C_m h^m,C_k h^k]\right\}, \, \mbox{ if }\, h<h^*, \label{B}
\end{eqnarray}
where $h^*$ is defined by:
\begin{equation}\label{h*}
\D h^* \equiv\left( \frac{C_k}{C_m}\right)^{\frac{1}{m-k}}.
\end{equation}
Then we showed in \cite{CMAM_2019} the two following results:
\begin{lemma}\label{Two_Steps}
Let $A$ and $B$ be the two events defined by (\ref{A}) and (\ref{B}) and let us assume they are independent. Then, the probability law of the event $\D\left\{X^{(m)}(h) \leq X^{(k)}(h)\right\}$ is given by:
\vspace{-0.2cm}
\begin{equation}\label{Heaviside_Prob}
\D Prob\left\{ X^{(m)}(h) \leq X^{(k)}(h)\right\} = \left |
\begin{array}{ll}
\hs 1 & \mbox{ if } \hs 0 < h < h^*, \medskip \\
\hs 0 & \mbox{ if } \hs h> h^*.
\end{array}
\right.
\end{equation}
\end{lemma}
Furthermore, if we replace the hypothesis of independency between $A$ and $B$ by considering the two random variables $X^{(i)}, (i=k \mbox{ or } i=m),$ independent and uniformly distributed $[0, C_i h^i]$, we also proved in \cite{CMAM_2019} the following theorem:
\begin{theorem}\label{The_nonlinear_law}
Let $u$ be the solution to the second order variational elliptic problem $\textbf{(VP}\textbf{)}$ defined in (\ref{VP}) and $u^{(i)}_h, (i=k \mbox{ or } i=m, k<m)$, the two corresponding Lagrange finite element $P_i$ approximations, solution to the approximated formulation $\textbf{(VP}\textbf{)}_{h}$ defined by (\ref{VP_h}).\\[0.1cm]
We assume the two corresponding random variables $X^{(i)}(h), (i=k \mbox{ or } i=m),$ are independent and uniformly distributed on $[0, C_i h^i]$, where $C_i$ are defined by (\ref{Ck_Cm}). \\[0.1cm]
Then, the probability of the event $\D\left\{ X^{(m)}(h) \leq X^{(k)}(h)\right\}$ is given by:
\begin{equation}\label{Nonlinear_Prob}
\D Prob\left\{ X^{(m)}(h) \leq X^{(k)}(h)\right\} = \left |
\begin{array}{ll}
\D \hs 1 - \frac{1}{2}\!\left(\!\frac{\!\!h}{h^*}\!\right)^{\!\!m-k} & \mbox{ if } \hs 0 < h \leq h^*, \\[0.5cm]
\D \hs \frac{1}{2}\!\left(\!\frac{h^*}{\!\!h}\!\right)^{\!\!m-k} & \mbox{ if } \hs h \geq h^*.
\end{array}
\right.
\end{equation}
\end{theorem}
The global shapes of the probabilistic laws (\ref{Heaviside_Prob}) and (\ref{Nonlinear_Prob}) are plotted in Figure \ref{Sigmoid} and new features of the relative finite elements accuracy are described in \cite{CMAM_2019}. Amongst other, these laws clearly indicate that there exist cases, (if $h>h*$ then $Prob\left\{ X^{(m)}(h) \leq X^{(k)}(h)\right\}\leq 0.5$), where $P_{m}$ finite elements \emph{probably} must be overqualified and a significant reduction of implementation time and execution cost can be obtained without loss of accuracy by implementing $P_k$ finite element.\sa
Another interesting property concerns the "sigmoid" probability law (\ref{Nonlinear_Prob}) and its relationship with the two steps law (\ref{Heaviside_Prob}). To prove this relationship, as a first step, we need some new features of the Lagrange finite element $P_k$. This is the purpose of the next section.
\section{$P_k$ canonical basis estimates}\label{propPk}
\noindent In this section we follow the definitions and properties of the $P_k$ finite element in $\R^n$ described by P. A. Raviart and J. M. Thomas in \cite{RaTho82}. \sa
Let us consider $K \subset \R^n$ a n-simplex which belongs to a regular mesh ${\mathcal T}_h$. Since a complete polynomial of order $k$ which belongs to $P_k(K)$ contains
\begin{equation}\label{Dim_Pk}
\D N \equiv \left(
  \begin{array}{c}
n+k \\
n
  \end{array}
\right)
= \frac{(n+k)!}{n!\,k!}
\end{equation}
terms, each n-simplex of the mesh ${\mathcal T}_h$ must be associated to $N$ independent degrees of freedom to assure the unisolvence of the finite element. \sa
It is convenient to carry out all analysis of n-simplexes in terms of the so-called barycentric coordinates $\lambda_1,\dots,\lambda_{n+1}$ which satisfy $\D \sum_{i=1}^{n+1}\lambda_i=1$.\sa
A regularly spaced set of points $M_{i_1,\dots,i_{n+1}}$ may be defined in a n-simplex by the barycentric coordinates values, namely:
\begin{equation}\label{Points_de_K}
\D M_{i_1,\dots,i_{n+1}} = \left(\frac{i_1}{k},\dots,\frac{i_{n+1}}{k}\right), \hs 0 \leq i_1,\dots,i_{n+1} \leq k,
\end{equation}
satisfying:
\begin{equation}\label{sum_ij_egal_k}
i_1 + \dots +i_{n+1}=k.
\end{equation}
One can verify that the number of points defined by (\ref{Points_de_K})-(\ref{sum_ij_egal_k}) is equal to $N$, the dimension of $P_k(K)$ in (\ref{Dim_Pk}). \sa
Therefore, we introduce the canonical basis of functions $p_{i_1,\dots,i_{n+1}}$ of the variables $(\lambda_1, \dots, \lambda_{n+1})$ which belongs to $P_k(K)$ defined by:
\begin{equation}\label{shape_function}
\D p_{i_1, \dots, i_{n+1}}(\lambda_1, \dots, \lambda_{n+1}) \equiv \prod_{j=1}^{n+1}P_{i_j}(\lambda_j),
\end{equation}
where the auxiliary polynomial $P_{i_j}(\lambda_j)$ is given by:
\begin{equation}\label{P_ij}
\D P_{i_j}(\lambda_j) \equiv \left |
\begin{array}{ll}
\D \hs \prod_{c_j=1}^{i_j}\left(\frac{k \lambda_j - c_j +1}{c_j}\right), & \mbox{ if } \hs i_j \geq 1, \vspace{0.1cm} \\
\D\hs 1, & \mbox{ if } \hs i_j = 0.
\end{array}
\right.
\end{equation}
$P_{i_j}$ is clearly a polynomial of order $i_j$ in $\lambda_j$, and therefore, due to condition (\ref{sum_ij_egal_k}), $p_{i_1, \dots, i_{n+1}}$ given by (\ref{shape_function}) is a polynomial of order $k$. \sa
In the sequel, notice that we will also use a simple index numbering to substitute the multi-index numbering. It will be the case for the $N$ points $M_{i_1,\dots,i_{n+1}}$ simply denoted $(M_i)_{i=1,N}$, as well as for the $N$ canonical functions $p_{i_1,\dots,i_{n+1}}$ denoted $(p_i)_{i=1,N}$, and so on. \sa
Therefore, the main property of the canonical basis of functions $p_i \equiv p_{i_1,\dots,i_{n+1}}$ defined in (\ref{shape_function}) is that for a given set of $N$ values $\varphi_i \equiv \varphi_{{i_1, \dots, i_{n+1}}}$ known at the $N$ points $M_i \equiv M_{i_1, \dots, i_{n+1}}$, the polynomial $Q$ in $P_k(K)$ given by:
\begin{eqnarray}
\D\forall M \in K: Q(M) & = & Q(\lambda_1, \dots, \lambda_{n+1}) \nonumber \\[0.2cm]
& = & \sum_{i_1 + \dots + i_{n+1}=k} \!\!\!\! \varphi_{i_1, \dots, i_{n+1}} p_{i_1, \dots, i_{n+1}} \! (\lambda_1, \dots, \lambda_{n+1}) \nonumber \\[0.2cm]
& = & \hspace{0.7cm}\sum_{i=1}^{N} \varphi_i p_i(\lambda_1, \dots, \lambda_{n+1}),
\end{eqnarray}
is the unique one in $P_k(K)$ such that $Q(M_i)= \varphi_i$. \sa
The following result concerns the features of the canonical basis $(p_i)_{i=1,N}$ defined by (\ref{shape_function}), where $N$ is given by (\ref{Dim_Pk}), regarding the semi-norm $|.|_{m,p,K}$ in $W^{m,p}(K)$ in the particular cases $p=2$ and $m=0$ or $m=1$.\sa
First of all, we remark that the structure of the elementary polynomials $P_{ij}$ defined by (\ref{P_ij}) looks like the numerator of the famous Lagrange polynomials. Then, we will establish the first estimate.
\begin{lemma}\label{Numerator_Lagrange}
Let $[a,b], (a<b),$ be a given interval and $N_p$ a given non-zero integer. We consider a set of $N_p+1$ uniform distributed points $x_j, (j=0,\dots,N_p),$ in $[a,b]$ defined by:
\begin{equation}
\D \forall j=0,\dots,N_p: x_j=a+jh, h=\frac{b-a}{N_p}.
\end{equation}
Let also $\Pi$ be the function defined on $[a,b]$ by: $\D\Pi(x)= \prod_{j=0}^{N_p}(x-x_j)$. \sa
Then,
\begin{equation}\label{Pi_control}
\forall x \in [a,b]: |\Pi(x)| \leq (N_p+1)!\,h^{N_p+1}.
\end{equation}
\end{lemma}
\begin{prooff} $\frac{}{}$ \sa
First, remark that for $x=b$, $\Pi(b)$ vanishes, so that (\ref{Pi_control}) is satisfied. Hence, let $x$ be a fixed value in the interval $[a,b[$. It exists a unique $i\in\{0,\dots,N_p-1\}$ such that $x\in[x_i,x_{i+1}[$. Therefore, we write the function $\Pi(x)$ as follows:
\begin{equation} \label{Pi_with_i}
\D \Pi(x) = (x-x_0)\dots(x-x_i)(x-x_{i+1})\dots(x-x_{N_p}).
\end{equation}
Furthermore, we have the following inequalities:
\begin{eqnarray}
& & \D\forall j=0,\dots,i : |x-x_j| \leq (i+1-j)h, \label{Majo_1}\\[0.2cm]
& & \D\forall j=1,\dots,N_p-i : |x-x_ {i+j}| \leq jh \leq (i+j+1)h.\label{Majo_2}
\end{eqnarray}
Therefore, $\Pi(x)$ written in (\ref{Pi_with_i}) can be controlled by the help of (\ref{Majo_1}) and (\ref{Majo_2}) by:
\begin{equation}\label{Control_Pi}
\D|\Pi(x)| \leq \left[(i+1)h\times(ih)\dots (2h)\times (h)\right]\times\left[(i+2)h\times(i+3)h\dots(N_p+1)h\right],
\end{equation}
and after reorganizing the right side of (\ref{Control_Pi}), we get (\ref{Pi_control}).
\end{prooff}
The following lemma gives us the first point-to-point estimates for the polynomials $p_i$ defined by (\ref{shape_function}).
\begin{lemma}
Let $p_i, (i=1,\dots,N)$, be the basis functions of the space of polynomials $P_k(K)$ which are defined by (\ref{shape_function})-(\ref{P_ij}). \sa
Then, $\forall i=1,\dots,N, \,\forall l=1,\dots,n+1:$
\begin{equation}\label{Local_estimate_pi_and_Dpi}
|p_i(\lambda_1, \dots, \lambda_{n+1})|\leq k^{n+1}, \hs \left|\frac{\partial p_i}{\partial\lambda_l}(\lambda_1, \dots, \lambda_{n+1})\right|\leq k^{n+2}.
\end{equation}
\end{lemma}
\begin{prooff} $\frac{}{}$
Let us introduce the integer $n_i$, $(0 \leq n_i\!\leq n+1)$, which corresponds to the number of polynomials $P_{i_j}(\lambda_j)$ such that:
\vspace{-0.2cm}
\begin{eqnarray}
\forall j=1,\dots,n_i, \,(n_i \geq 1),\,& : & P_{i_j}(\lambda_j)=P_{1}(\lambda_j)=k\lambda_j, (i_j=1), \label{Num_1} \\[0.1cm]
\forall j=n_i+1,\dots,n+1, \,(n_i \leq n)\, & : & P_{i_j}(\lambda_j)=\frac{k\lambda_j(k\lambda_j - 1)\dots(k\lambda_j - i_j+1)}{i_j!}, (i_j>1). \label{Num_2}
\end{eqnarray}
When $n_i=0$, then $p_i$ has the following structure:
\begin{equation}\label{pi_0}
\D p_i(\lambda_1,\dots,\lambda_{n+1}) = \prod_{j=1}^{n+1}\left(\frac{k\lambda_j(k\lambda_j - 1)\dots(k\lambda_j - i_j+1)}{i_j!}\right), (i_j >1, \forall j=1,\dots,n+1),
\end{equation}
and when $n_i=n+1$, then $p_i$ corresponds to:
\begin{equation}\label{pi_n_i}
\D p_i(\lambda_1,\dots,\lambda_{n+1}) = \prod_{j=1}^{n+1}(k\lambda_{i_j}).
\end{equation}
$\bullet$ Let us begin by fixing a given value of $n_i, (1 \leq n_i \leq n)$.\sa
Concerning the control of the polynomials $p_i$, we split it into two groups of elementary polynomials $P_{i_j}$ as follows:
\begin{equation}\label{pi_1}
\D p_i(\lambda_1,\dots,\lambda_{n+1}) = \prod_{j=1}^{n_i}(k\lambda_{i_j}).\!\!\!\prod_{j=n_i+1}^{n+1}\left(\frac{k\lambda_j(k\lambda_j - 1)\dots(k\lambda_j - i_j+1)}{i_j!}\right).
\end{equation}
Now, on the first hand, the barycentric functions $\lambda_j, (j=1\dots,n+1)$ satisfy: $$\forall M\in K: 0\leq\lambda_j(M)\leq 1.$$
On the other hand, by applying lemma \ref{Numerator_Lagrange} by setting $x=k\lambda_j$ and $h=1$, we have the following estimate:
\begin{equation}\label{pi_2}
\D \forall j=n_i+1,\dots,n+1: |P_{i_j}(\lambda_j)| = \left|\frac{k\lambda_j(k\lambda_j - 1)\dots(k\lambda_j - i_j+1)}{i_j!}\right| \leq  1,
\end{equation}
and finally,
\begin{equation}\label{pi_3}
\D |p_i(\lambda_1,\dots,\lambda_{n+1})| \leq \left|\prod_{j=1}^{n_i}(k\lambda_{i_j})\right| \leq k^{n_i} \leq k^{n},
\end{equation}
as $n_i \leq n$.\sa
Let us now consider the partial derivative $\D\frac{\partial p_i}{\partial\lambda_l}$, for a given pair of non zero integers $(i,l)$. \sa
By (\ref{shape_function}) we can write the concerned partial derivative as:
\begin{equation}\label{Dpi_il_on_Dlambda_l_1}
\D\frac{\partial p_i}{\partial\lambda_l} =P_{i_1}\dots\frac{\partial P_{i_l}}{\partial\lambda_l} \dots P_{i_{n+1}}.
\end{equation}
Thus, two cases have to be considered. The first one corresponds to the case when $P_{i_l}$ is a single monomial $(i_l=1)$: $P_{i_l}(\lambda_l)= k\lambda_l$.\sa
Therefore, (\ref{Dpi_il_on_Dlambda_l_1}) gives:
\begin{equation}\label{Dpi_il_on_Dlambda_l_2}
\D\frac{\partial p_i}{\partial\lambda_l} =k\,P_{i_1}\dots P_{i_{l}-1}P_{i_{l}+1}\dots P_{i_{n+1}},
\end{equation}
and similarly to (\ref{pi_1})-(\ref{pi_3}), we get the following estimate:
\begin{equation}\label{Dpi_il_on_Dlambda_l_3}
\D\left|\frac{\partial p_i}{\partial\lambda_l}\right| \leq k.k^{n_i-1}\leq k^{n}.
\end{equation}
Let us now consider the case when $i_l>1$. It means that the polynomial $P_{i_l}$ has the structure of (\ref{pi_1}), composed at least by two monomials. Then, its partial derivative with respect to $\lambda_l$ is equal to:
\begin{equation}\label{Dpi_il_on_Dlambda_l_4}
\D\frac{\partial P_{i_l}}{\partial\lambda_l} =\frac{1}{i_l !}\left[\frac{}{}\!\!k(k\lambda_l-1)\dots(k\lambda_l-i_l+1) + \dots + k\lambda_l(k\lambda_l-1)\dots(k\lambda_l-i_l+2)k\right].
\end{equation}
So, by using the same arguments we implemented to upper bound the function $\Pi$ defined in lemma \ref{Numerator_Lagrange}, on can increase each term of the right hand side of (\ref{Dpi_il_on_Dlambda_l_4}) to finally obtain:
\begin{equation}\label{Dpi_il_on_Dlambda_l_5}
\D\left|\frac{\partial P_{i_l}}{\partial\lambda_l}\right| \leq \frac{k}{i_l !}\left[\frac{}{}\!\!\, i_l!\,i_l \right]\leq k^2,
\end{equation}
as $\forall l=1$ to $n+1$: $i_l \leq k$.\sa
Finally we get the estimate for the partial derivative of $p_i$ with respect to $\lambda_l$:
\begin{equation}\label{Dpi_il_on_Dlambda_l_6}
\D\left|\frac{\partial p_i}{\partial\lambda_l}\right| \leq k^2 \, k^{n_i} \leq k^{n+2}.
\end{equation}
Let us now consider the two cases when $n_i=0$ or $n_i=n+1$.\sa
$\bullet$ If $n_i=0$ then, due to (\ref{pi_0}) we have:
\begin{equation}\label{pi_00}
\D |p_i(\lambda_1,\dots,\lambda_{n+1})| \leq  1 \leq \ k^{n},
\end{equation}
thanks to (\ref{pi_2}).\sa
In the same way, we have the following inequalities:
\begin{equation}\label{Dpi_il_on_Dlambda_ni=0}
\D\left|\frac{\partial p_i}{\partial\lambda_l}\right| \leq k^2 \leq k^{n+2},
\end{equation}
where we used (\ref{Dpi_il_on_Dlambda_l_4})-(\ref{Dpi_il_on_Dlambda_l_5}) which correspond to the present case.\sa
$\bullet$ If $n_i=n+1$ then $p_i$ is given by (\ref{pi_n_i}) and we have:
\begin{equation}
\D |p_i(\lambda_1,\dots,\lambda_{n+1})| \leq k^{n+1} \mbox{ and }  \left|\frac{\partial p_i}{\partial\lambda_l}\right| \leq k^{n+1},
\end{equation}
due to the basic barycentric functions features.\sa
Therefore, from all the above upper bounds, we get (\ref{Local_estimate_pi_and_Dpi}).
\end{prooff}
We can now get two estimates of the canonical basis $(p_i)_{i=1,N}$ with respect to the semi-norms $|p_i|_{0,2,K}$ and $|p_i|_{1,2,K}$.
\begin{lemma}\label{Estimation_pi}
Let $(p_i)_{i=1,N}$ be the canonical basis defined in (\ref{shape_function}). Then, if $\D k >\frac{n}{2},$ we have:
\begin{equation}\label{Norm_0_2_and_Norm_1_2}
\D \sum_{i=1}^{N}|p_i|_{0,2,K} = {\mathcal O\!}\left(\frac{}{}\!\!k^{n+1} (k+n)^{n}\right) \mbox{ and } \sum_{i=1}^{N}|p_i|_{1,2,K} = \D\frac{\D{\mathcal O}\left(\frac{}{}\!\!k^{n+2} (k+n)^{n}\right)}{\rho_{K}},
\end{equation}
where ${\mathcal O}$ denotes Landau's notation and $\rho_{K}$ the diameter of the largest inscribed sphere within $K$.\\
\end{lemma}
\begin{prooff} $\frac{}{}$ \sa
\noindent \hspace{-0.1cm}$\blacktriangleright$ Let us begin with the estimate of $p_i$ with respect to the semi-norm $|.|_{0,2,K}$. \sa
From the local estimate of $p_i$ given by (\ref{Local_estimate_pi_and_Dpi}), we directly get the $|.|_{0,2,K}-$semi-norm for each polynomial $p_i, (1 \leq i\leq N),$ as follows:
\begin{equation}\label{Norm_L2_pi_00_1}
|p_i|_{0,2,K} \leq \sqrt{\mbox{mes}(K)} k^{n+1},
\end{equation}
and by aggregating on all the $N$ basis function $p_i$, we get:
\begin{eqnarray}
\D \sum_{i=1}^{N}|p_i|_{0,2,K} & \leq & \sqrt{\mbox{mes}(K)} \,\, k^{n+1}\, \frac{(n+k)!}{n!\,k!}, \nonumber\\[0.2cm]
\D & \leq & \sqrt{\mbox{mes}(K)}\,(k+n)^n k^{n+1}, \label{Norm_L2_pi_0_2}
\end{eqnarray}
as $n\geq 1$ and where we used the value of $N$ corresponding to the dimension of the space $P_k(K)$ given by (\ref{Dim_Pk}). \sa
Finally, with (\ref{Norm_L2_pi_0_2}) we get the first estimate of (\ref{Norm_0_2_and_Norm_1_2}). \sa
\noindent \hspace{-0.1cm}$\blacktriangleright$ Let us prove now the second estimate of (\ref{Norm_0_2_and_Norm_1_2}) with respect to the semi-norm $|.|_{1,2,K}$. \sa
Due to remark 2.2 in R. Arcangeli and J. L. Gout \cite{Arcangeli_Gout}, for each canonical basis function $p_i$, if $\D k > \frac{n}{2},$ we have:
\begin{equation}\label{Norm_Grad_L2_pi_0}
\D |p_i|_{1,2,K} \leq \frac{1}{\rho_{K}}\left\{\int_{K}\left[\sum_{j=1}^{n}\left|\frac{\partial p_i}{\partial x^j}(x)\right|\right]^{2}dx\right\}^{\frac{1}{2}},
\end{equation}
where $\rho_{K}$ is the supremum of the diameters of the inscribed spheres within the n-simplex $K$.\sa
Moreover, each partial derivative $\D \frac{\partial p_i}{\partial x^j}$ can be computed using the chain rule as follows:
\begin{equation}
\D \frac{\partial p_i}{\partial x^j} = \sum_{l=1}^{n+1}\frac{\partial p_i}{\partial \lambda_l}\frac{\partial \lambda_l}{\partial x^j},
\end{equation}
where each partial derivative $\D \frac{\partial \lambda_l}{\partial x^j}$ is a constant $\Lambda^{(l)}_j$ that does not depend on $k$,  $\lambda_l$ being a polynomial of degree at most equal to one. \sa
So, we have:
\begin{equation}
\D \frac{\partial p_i}{\partial x^j} = \sum_{l=1}^{n+1}\Lambda^{(l)}_j\frac{\partial p_i}{\partial \lambda_l}.
\end{equation}
Consequently, from (\ref{Norm_Grad_L2_pi_0}), we get:
\begin{eqnarray}
\D |p_i|^2_{1,2,K} & \leq & \frac{1}{\rho^2_{K}}\int_{K}\left[\sum_{j=1}^{n}\left|\sum_{l=1}^{n+1}\Lambda^{(l)}_j\frac{\partial p_i}{\partial \lambda_l}\right|\right]^{2}dx
\\[0.2cm]
& \leq & \left(\frac{\D n\Lambda}{\rho_{K}}\right)^2\int_{K} \left|\sum_{l=1}^{n+1}\frac{\partial p_i}{\partial \lambda_l}\right|^{2}dx, \label{Norm_Grad_L2_pi_1} \\[0.2cm]
& \leq & \mbox{mes}(K)\left(\frac{\D n(n+1)\Lambda \, k^{n+2}}{\rho_{K}}\right)^2, \label{Norm_Grad_L2_pi_2}
\end{eqnarray}
where we set $\D \Lambda\equiv\!\!\!\max_{1 \leq j \leq n \atop \vspace{0.3cm} \hspace{0.3cm} 1 \leq l \leq n+1} \!\!\left|\Lambda^{(l)}_j\right|$ and due to (\ref{Dpi_il_on_Dlambda_l_6}).\sa
By aggregating (\ref{Norm_Grad_L2_pi_2}) on the $N$ basis functions $p_i$ we finally get:
\begin{equation}\label{Norm_Grad_L2_pi_3}
\D \sum_{i=1}^{N}|p_i|_{1,2,K} \leq \sqrt{\mbox{mes}(K)}\,\frac{\D n(n+1)\Lambda}{\rho_{K}}\, (k+n)^n k^{n+2},
\end{equation}
which corresponds to the second estimate of (\ref{Norm_0_2_and_Norm_1_2}).
\end{prooff}
\begin{remark}
We notice that for applications the condition $\D k > \frac{n}{2}$ holds for the dimension $n=1$ when $k\geq 1$, but if $n=2$ or $n=3$ this requires $k\geq 2$. Consequently, the case of the finite element $P_1$ could be still considered by using other results of \cite{Arcangeli_Gout}, as we will mention later, (see Theorem \ref{C_k_h*q}).
\end{remark}
\noindent The two estimates (\ref{Norm_0_2_and_Norm_1_2}) will now be used to determine the asymptotic behavior of the probability distribution (\ref{Nonlinear_Prob}) of theorem~\ref{The_nonlinear_law}.
\section{Asymptotic limit of the "sigmoid" probability distribution}\label{Asymptotic_limit}
\begin{figure}[h]
  \centering
  \includegraphics[width=12cm]{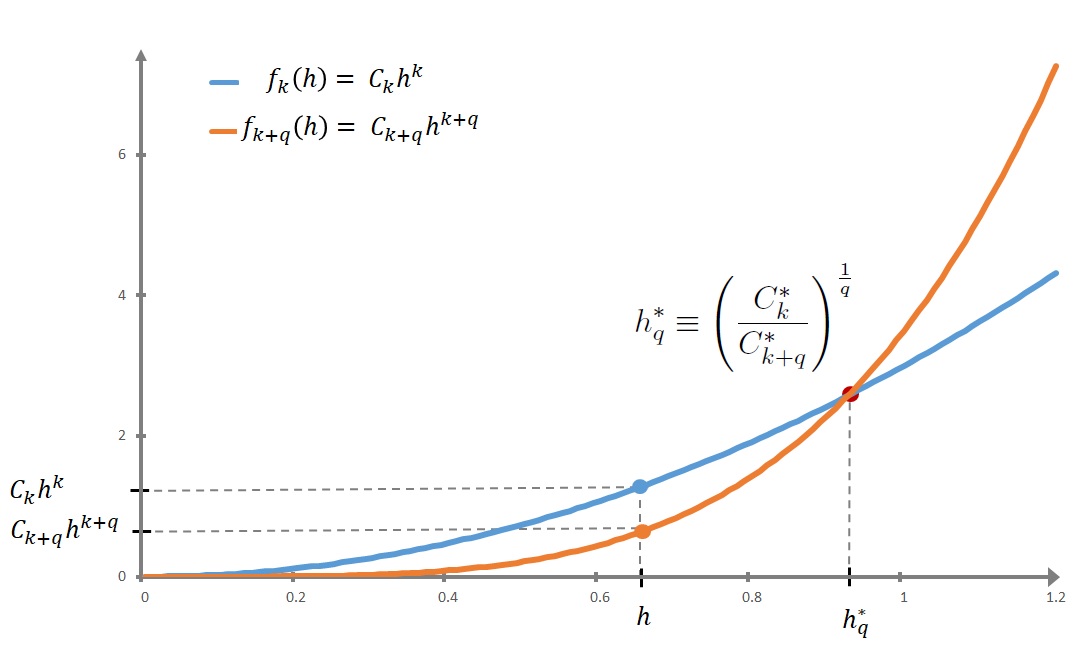}
  \caption{Relative positions between $f_k$ and $f_{k+q}$ curves.} \label{Xk+q_vs_Xk}
\end{figure}

\noindent As we already mentioned, the probability distribution (\ref{Nonlinear_Prob}) has been approximated by the stepwise law (\ref{Heaviside_Prob}) if one assumes the independency between the events $A$ and $B$ defined by (\ref{A}) and (\ref{B}). \sa
Conversely, here we will study the behavior of the non linear law (\ref{Nonlinear_Prob}) when $q \equiv m-k$ goes to infinity. This study is not only theoretical. It is clearly related to the well-known question, namely,  in which way high order finite element methods can solve partial differential equations more efficiently than low order methods. More precisely, how large of a polynomial degree is beneficial? Here again, we have chosen to treat the problem {\em via} a probabilistic approach, handling the uncertainties (randomness of the data, of the mesh, etc.) by random variables. Note that often in the applications, one considers cases when $k=1$ or $2$, whereas the high order degree $m$ is around $20-25$, see for instance \cite{Mitc15}. \sa
More precisely, let us give a fixed value of $k$. \sa
Then, we define the sequence of functions $\D\left(\mathcal{P}_{q}(h)\right)_{q \in \N^\star}$ by:
\begin{equation}\label{Suite_Pn}
\D \mathcal{P}_{q}(h) \equiv Prob\left\{ X^{(k+q)}(h) \leq X^{(k)}(h)\right\},
\end{equation}
where $\D Prob\left\{ X^{(k+q)}(h) \leq X^{(k)}(h)\right\}$, following (\ref{Nonlinear_Prob}), is given by:
\begin{equation}\label{Nonlinear_Prob_Truncate}
\D Prob\left\{ X^{(k+q)}(h) \leq X^{(k)}(h)\right\} = \left |
\begin{array}{ll}
\D \hs 1 - \frac{1}{2}\!\left(\!\frac{\!\!h}{h_q^*}\right)^{\!\!q} & \D\mbox{ if } \hs  0 < \frac{\!\!h}{h_q^*} \le 1, \medskip \\
\D \hs \frac{1}{2}\!\left(\!\frac{h_q^*}{\!\!h}\!\right)^{\!\!q} & \D\mbox{ if } \hs  \frac{\!\!h}{h_q^*} \ge 1,
\end{array}
\right.
\end{equation}
and where $(h^*_{q})_{q\in\N^\star}$ is the sequence defined by:
\begin{equation}\label{h*n}
\D h^*_{q} \equiv\left(\frac{C_k}{C_{k+q}}\right)^{\frac{1}{q}}.
\end{equation}
As one can see the critical value $h^*_{q}$ strongly depends on $q$, among others, by the constant $C_{k+q}$.\sa
To this end, we will firstly determine an estimate of the constant $\mathscr{C}_k$ defined by (\ref{Constante_01}), relatively to the finite element $P_k$.
This is the purpose of the following theorem:
\begin{theorem}\label{C_k_h*q}
For $\D k > \frac{n}{2}$, let $\mathscr{C}_k$ be the constant introduced in the error estimate (\ref{Constante_01}). Then, the following estimation holds: %
\begin{equation}\label{C_k_Estimation}
\D  \mathscr{C}_k = {\mathcal O}\left(\frac{(k+n)^n k^{\,n+2}}{(k-1)!\,\left(\!k-\D\frac{n}{2}\right)}\right).
\end{equation}
\end{theorem}
\begin{prooff}
The proof of this theorem is based on the paper of R. Arcangeli and J.L. Gout \cite{Arcangeli_Gout},  itself an extension of the one of P.G. Ciarlet and P.A. Raviart \cite{Ciarlet_Raviart}. \sa
To this end, let us firstly recall the conditions of theorem 2.1 of R. Arcangeli and J.L. Gout. \sa
Let $\Omega$ be an open bounded and non empty convex subset of $\R^n$ and $\Gamma$ its Lipschitz boundary. We assume that $\Sigma=\{a_i\}_{i=1,N}$ is a $P-$unisolvent set of points which belong to $\bar{\Omega}$, where $P$ denotes a space of finite dimension such that $P_k \subset P \subset C^{k}(\bar{\Omega})$, and $P_k$ the space of polynomial functions of degree less than or equal to $k$. \sa
Then, for all $u\in W^{k+1,p}(\Omega)$ and for all integer $\nu \geq 0$ such that
\begin{equation}\label{cond_parametrique}
\D k+1 > \nu + \frac{n}{p},
\end{equation}
we have:
\begin{eqnarray}
\D |u -\Pi_h u|_{\nu,p,\Omega} & \leq & \frac{1}{(k-\nu)!\left(k+1-\nu-\D\frac{n}{p}\right)}|u|_{k+1,p,\Omega}\,h^{k+1-\nu}  \nonumber \\[0.3cm]
& + & \D\frac{\left(\D\sum_{i=1}^{N}|p_i|_{\nu,p,\Omega} \!\right)}{[\mbox{mes}(\Omega)]^{1/p}\,k!\!\left(k+1-\D\frac{n}{p}\right)}|u|_{k+1,p,\Omega}\,h^{k+1}, \label{Estimation_Arc_Gout}
\end{eqnarray}
where $\D |.|_{\nu,p,\Omega}$ denotes the usual semi-norm in the Sobolev spaces $W^{\nu,p}(\Omega)$, $\Pi_h$ the classical Lagrange interpolation which consists to interpolate the set of points $\Sigma$ in $\R^n$ by a polynomial function of a given degree $k$, and $(p_i)_{i=1,N}$ the unique functions such that
$$p_{i}(Mj)=\delta_{ij}, \forall\, M_j \in \Sigma, \forall \, 1\leq i, j\leq N,$$
where $\delta_{ij}$ denotes the Kr\"onecker's symbol. \sa
Here, for our objectives, we write (\ref{Estimation_Arc_Gout}) in the particular following conditions:
\begin{itemize}
\item $\Omega=K_\mu$, where $K_\mu, (1 \leq \mu \leq N_K),$ is a n-simplex which belongs to a given regular mesh ${\mathcal T}_h$.
\item $u$ is the exact solution to the variational formulation \textbf{(VP)} defined in (\ref{VP}).
\item The set of points $\Sigma$ in $\R^n$ corresponds to the $P_k$ finite element degrees of freedom defined on the n-simplex $K_{\mu}$ defined by (\ref{Points_de_K}) and (\ref{sum_ij_egal_k}).
\item The interpolation operator $\Pi_h$ is replaced by $\Pi_{K_{\mu}}$, the local Lagrange interpolation operator.
\end{itemize}
Then, we choose in (\ref{Estimation_Arc_Gout}) $p=2$, $\nu=0$ and $\nu=1$ which implies that $k > \D \frac{n}{2}$ due to (\ref{cond_parametrique}), or equivalently, $k \geq 2$ for a problem set in dimension $n\geq 2$.\sa
The case of the finite element $P_1$ in dimension $n\geq 2$ could also be considered by adapting our theorem with another result from R. Arcangeli and J.L. Gout (see remark 2.3 and theorem 1.1 in \cite{Arcangeli_Gout}). \sa
\noindent So, we get the following inequalities:\sa
$\blacktriangleright$ For $\nu=0$ we have:
\begin{eqnarray}
\D \hspace{-7.5cm}\forall\, K_{\mu} \in {\mathcal T}_h, 1 \leq \mu \leq N_K &\hspace{-0.1cm} : \hspace{-0.1cm}& \nonumber
\end{eqnarray}
\vspace{-0.8cm}
\begin{eqnarray}
\D |u -\Pi_{K_{\mu}}u|_{0,K_{\mu}} & \leq & \frac{1}{k!\left(k+1-\D\frac{n}{2}\right)}|u|_{k+1,K_{\mu}}\,h_{K_{\mu}}^{k+1}  \nonumber \\[0.3cm]
& + & \frac{\left(\D\sum_{i=1}^{N}|p_i|_{0,K_{\mu}} \right)}{\left[\mbox{mes}(K_{\mu})\right]^{1/2} k!\left(k+1-\D\frac{n}{2}\right)}|u|_{k+1,K_{\mu}}\,h_{K_{\mu}}^{k+1}, \label{Estimation_Arc_Gout_m=0_0}
\end{eqnarray}
which becomes:
\begin{equation}\label{Estimation_Arc_Gout_m=0}
\D |u -\Pi_{K_{\mu}}u|_{0,K_{\mu}} \leq  \left[\frac{1 + (k+n)^n k^{n+1}}{k!\left(k+1-\D\frac{n}{2}\right)}\right]|u|_{k+1,K_{\mu}}\,h_{K_{\mu}}^{k+1},
\end{equation}
due to (\ref{Norm_L2_pi_0_2}). \\[0.2cm]
$\blacktriangleright$ In the same way, for $\nu=1$, we have:
\begin{eqnarray}
\D |u -\Pi_{K_{\mu}}u|_{1,K_{\mu}} & \leq & \frac{1}{(k-1)!\left(k-\D\frac{n}{2}\right)}|u|_{k+1,K_{\mu}}\,h_{K_{\mu}}^{k} \nonumber \\[0.3cm]
\D & + & \frac{\D\left(\sum_{i=1}^{N}|p_i|_{1,K_{\mu}} \right)}{\left[\mbox{mes}(K_{\mu})\right]^{1/2}k!\left(k+1-\D\frac{n}{2}\right)}|u|_{k+1,K_{\mu}}\,h_{K_{\mu}}^{k+1}. \label{Estimation_Arc_Gout_m=1}
\end{eqnarray}
which leads to:
\begin{eqnarray}
\D |u -\Pi_{K_{\mu}}u|_{1,K_{\mu}} & \leq & \frac{1}{(k-1)!}\frac{1}{\left(k-\D\frac{n}{2}\right)}|u|_{k+1,K_{\mu}}\,h_{K_{\mu}}^{k} \nonumber \\[0.2cm]
& + &\frac{n(n+1)\Lambda}{\rho_{K_{\mu}}} \frac{(k+n)^n k^{n+2} }{k! \left(k+1-\D\frac{n}{2}\right)}|u|_{k+1,K_{\mu}}\,h_{K_{\mu}}^{k+1},
\end{eqnarray}
due to (\ref{Norm_Grad_L2_pi_3}), and finally:
\begin{equation}\label{Estimation_Arc_Gout_m=1_1}
\D |u -\Pi_{K_{\mu}}u|_{1,K_{\mu}} \leq \left[\frac{1+ \sigma n(n+1) \Lambda(k+n)^n k^{n+2}}{(k-1)!\,\left(k-\D\frac{n}{2}\right)}\right] |u|_{k+1,K_{\mu}}\,h_{K_{\mu}}^{k},
\end{equation}
where we introduced the number $\sigma \geq 1$ such that $\D\frac{h_{K_{\mu}}}{\rho_{K_{\mu}}}\leq \sigma$, $\forall\,K_{\mu} \in {\mathcal T}_h$.\sa
Therefore, by the help of (\ref{Estimation_Arc_Gout_m=0}) and (\ref{Estimation_Arc_Gout_m=1_1}), we get the following estimate of the local interpolation error with respect to the $H^1-$norm, using that the mesh ${\mathcal T}_h$ is regular, and by setting $\D h \equiv \max_{K_{\mu} \in {\mathcal T}_h} h_{K_{\mu}}$:
\begin{eqnarray}
\D \|u -\Pi_{K_{\mu}} u\|^2_{1,K_{\mu}} & \leq & \left[\frac{1 + (k+n)^n k^{n+1}}{k!\left(k+1-\D\frac{n}{2}\right)}\right]^2\!
|u|^2_{k+1,K_{\mu}}\,h_{K_{\mu}}^{2(k+1)} \nonumber \\[0.2cm]
& + & \left[\frac{1+ \sigma n(n+1) \Lambda(k+n)^n k^{n+2}}{(k-1)!\,\left(k-\D\frac{n}{2}\right)}\right]^2\! |u|^2_{k+1,K_{\mu}}\,h_{K_{\mu}}^{2k}. \nonumber
\end{eqnarray}
Then, we get:
\begin{eqnarray}
\D \|u -\Pi_{K_{\mu}} u\|_{1,K_{\mu}} & \leq & \left[C(\bar{\Omega},\sigma,\Lambda,n)\frac{(k+n)^n k^{n+2}}{(k-1)!\,\left(k-\D\frac{n}{2}\right)}\right]\! |u|_{k+1,K_{\mu}}\,h^{k}, \label{Lambda_000}
\end{eqnarray}
where we introduced the constant $C(\bar{\Omega},\sigma,\Lambda,n)$ defined by:
\begin{equation}\label{C}
\D C(\bar{\Omega},\sigma,\Lambda,n) \equiv 1+ 2\,\mbox{diam}(\bar{\Omega}) + \sigma n(n+1) \Lambda.
\end{equation}
Therefore, by the help of (\ref{Lambda_000}), we get for the whole domain $\Omega$ the following estimate of the interpolation error:
\begin{eqnarray}
\D \|u -\Pi_h u\|_{1,\Omega} & = & \left(\sum_{K_{\mu}\in {\mathcal T}_h}\|u -\Pi_{K_{\mu}}u\|_{1,K_{\mu}}^2\!\!\right)^{\!\!1\!/2} \nonumber \\
& \leq & C(\bar{\Omega},\sigma,\Lambda,n)\frac{(k+n)^n k^{n+2}}{(k-1)!\,\left(k-\D\frac{n}{2}\right)}\left(\!\sum_{K_{\mu}\in {\mathcal T}_h}|u|^2_{k+1,K_{\mu}}\!\!\right)^{\!\!1\!/2}\!\!\!\!\!\!h^{k},\nonumber\\[0.4cm]
& \leq & C(\bar{\Omega},\sigma,\Lambda,n)\frac{(k+n)^n k^{n+2}}{(k-1)!\,\left(k-\D\frac{n}{2}\right)}\,|u|_{k+1,\Omega}\,h^{k}. \label{Lambda_0000}
\end{eqnarray}
Then, inequality (\ref{Lambda_0000}) leads to the estimate (\ref{C_k_Estimation}) if one takes into account the estimate of C\'ea's lemma \cite{RaTho82}. Indeed, consider the $H^1-$norm to measure the difference between the exact solution $u$ to the variational problem (\ref{VP}) and its approximation $u_h$ solution to (\ref{VP_h}), we have:
\begin{equation}\label{Erreur_approx_Erreur_global_interpol}
\|u-u_h\|_{1,\Omega} \hs \leq \hs \frac{M}{\alpha}\,\|u-\Pi_h u\|_{1,\Omega},
\end{equation}
where $M$  is the continuity constant and $\alpha$ the ellipticity constant of the bilinear form $a(\cdot,\cdot)$. \sa
As a consequence, by the help of (\ref{Lambda_0000}) we obtain that the constant $\mathscr{C}_k$ in (\ref{Constante_01}) satisfies:
\begin{equation}\label{C_k_Asympt}
\D  \mathscr{C}_k \leq \frac{MC(\bar{\Omega},\sigma,\Lambda,n)}{\alpha}\,\frac{(k+n)^n k^{n+2}}{(k-1)!\,\left(k-\D\frac{n}{2}\right)},
\end{equation}
which corresponds to (\ref{C_k_Estimation}).
\end{prooff}
For the sequel, we introduce the constant $\mathscr{C}^*_k$ defined by:
\begin{equation}\label{C_k_Asympt_2}
\D \mathscr{C}_k^* \equiv \frac{MC(\bar{\Omega},\sigma,\Lambda,n)}{\alpha}\,\frac{(k+n)^n k^{n+2}}{(k-1)!\,\left(k-\D\frac{n}{2}\right)},
\end{equation}
and the corresponding $h^*_q$ defined in (\ref{h*n}), in which we substitute $C_k$ by the corresponding value of $C^*_k$, that is:
\begin{equation}\label{h*q_2}
\D h^*_{q} \equiv\left(\frac{C^*_k}{C^*_{k+q}}\right)^{\frac{1}{q}} = \left(\frac{\mathscr{C}^*_k |u|_{k+1,\Omega}}{\mathscr{C}^*_{k+q}|u|_{k+q+1,\Omega}}\right)^{\frac{1}{q}}.
\end{equation}
As we are interested in the asymptotic behavior of $h^*_{q}$ defined by (\ref{h*q_2}) when $q$ goes $+\infty$, we will assume that the solution $u$ to the variational problem \textbf{(VP)} belongs to $H^r(\Omega), (\forall\, r\in \N)$.\sa
We are now in position to propose an estimate of the sequence $(h^*_q)_{q\in\N^\star}$ defined by (\ref{h*q_2}), when $q$ goes to infinity and the corresponding asymptotic limit of the sequence of functions defined by (\ref{Suite_Pn}).
\begin{theorem}\label{Conv_Simple}
Let us assume that solution $u$ of problem \textbf{(VP)} belongs to $H^{r}(\Omega), (\forall\, r\in \N)$.
Let also $(h^*_{q})_{q\in\N^\star}$ be the sequence defined by (\ref{h*q_2}) and $\D\left(\frac{}{}\!\!\mathcal{P}_{q}(h)\right)_{q\in\N^\star}$ the corresponding sequence of functions defined by (\ref{Suite_Pn}).\sa
For a fixed value of $k$, $(k > \D\frac{n}{2})$, if
\begin{equation}\label{Cond_ Ratio_Semi_Norm}
\D \lim_{q\rightarrow +\infty}\frac{|u|_{k+q+2,\Omega}}{|u|_{k+q+1,\Omega}} = l, (l\in\R^{*}_{+}),
\end{equation}
then,
\begin{equation}\label{equiv_h*q_inf_1}
h^*_q \underset{q \to +\infty}{\sim} \frac{1}{e\,l}\,q, \, \mbox{ and } \,\lim_{q\rightarrow +\infty}h^*_q = +\infty.
\end{equation}
Moreover, the sequence of functions $\D\left(\frac{}{}\!\!\mathcal{P}_{q}(h)\right)_{q\in\N^\star}$converges pointwise when $q$ goes to $+\infty$ to the function $\mathcal{P}_0$ defined on $\R^*_+$ by:
\begin{equation}\label{Heaviside_Prob*}
\mathcal{P}_{0}(h)= \left |
\begin{array}{ll}
1 & 0 \leq h < +\infty \,,\\
\D\frac{1}{2} & h = +\infty .
\end{array}
\right.
\end{equation}
\end{theorem}
\noindent \begin{prooff}$\frac{}{}$\\[0.1cm]
$\blacktriangleright$ Let us replace the expression of $\mathscr{C}^*_k$ defined by (\ref{C_k_Asympt_2}) in $h^*_q$ given by (\ref{h*q_2}). Then, we have following asymptotic behavior:
\begin{equation}\label{C_k_0}
\D\left(h^*_q\right)^q \underset{q \to +\infty}{\sim} \frac{(k+n)^n k^{n+2}}{(k-1)!\,\left(k-\D\frac{n}{2}\right)} \, \frac{(q+k-1)!\left(q+k-\D\frac{n}{2}\right)}{\D\left(\!\!\frac{}{}q+k+n\right)^n \left(q+k\right)^{n+2}}\,.\frac{|u|_{k+1,\Omega}}{|u|_{k+q+1,\Omega}}.
\end{equation}
However, from Stirling's formula, when $q$ goes to $+\infty$, we can specify the equivalent of $h^{*}_q$  given by inequality (\ref{C_k_0}):
\begin{eqnarray}
\frac{(q+k-1)!\left(q+k-\D\frac{n}{2}\right)}{\D\left(q+k\right)^{n+2}\left(\!\!\frac{}{}q+k+n\right)^n} & \underset{q \to +\infty}{\sim} & \frac{\sqrt{2\pi (q+k-1)}\D\left(\frac{q+k-1}{e}\right)^{(q+k-1)}\!\!\!\!\left(q+k-\D\frac{n}{2}\right)}{\left(q+k\right)^{n+2}\left(\!\!\frac{}{}q+k+n\right)^n}, \nonumber \\[0.2cm]
& \underset{q \to +\infty}{\sim} & \frac{\sqrt{2 \pi}(q+k-1)^{(q+k-\frac{1}{2})}}{e^{q+k-1}}\frac{1}{\left(\!\!\frac{}{}q+k\right)^{2n+1}}, \nonumber \\[0.2cm]
& \underset{q \to +\infty}{\sim} & \sqrt{2 \pi}\,\frac{(q+k)^{q+k-2n-\frac{3}{2}}}{e^{q+k}}. \label{h*q_Asympt}
\end{eqnarray}
Then, (\ref{h*q_Asympt}) in (\ref{C_k_0}) leads to:
\begin{equation}\label{Truc}
\D\left(h^*_q\right)^q \underset{q \to +\infty}{\sim} \Theta \, e^{-(q+k)}(q+k)^{q+k-2n-\frac{3}{2}}\,.\frac{|u|_{k+1,\Omega}}{|u|_{k+q+1,\Omega}},
\end{equation}
where we introduced the constant $\Theta$ independent of $q$ defined by:
\begin{equation}
\Theta \equiv \sqrt{2 \pi}\frac{(k+n)^n k^{n+2}}{(k-1)!\,\left(k-\D\frac{n}{2}\right)}.
\end{equation}
Moreover, as we assume condition (\ref{Cond_ Ratio_Semi_Norm}), if we introduce the two sequences $(v_q)_{q\in\N}$ and $(w_q)_{q\in\N}$ as follows:
\begin{equation}
\forall\, q\in\N: v_q\equiv \ln\,|u|_{k+q+1,\Omega}, \hs w_q \equiv q,
\end{equation}
then, the ratio $r_q$ defined by:
\begin{equation}
\D r_q \equiv \frac{v_{q+1}-v_q}{w_{q+1}-w_q} = \ln\left(\frac{|u|_{k+q+2,\Omega}}{|u|_{k+q+1,\Omega}}\right),
\end{equation}
has a limit $L\equiv \ln\, l \in\R$, when $q$ goes to $+\infty$: $\D\lim_{q\rightarrow +\infty}r_q = L$.\sa
As a consequence, due to Stolz-Cesaro theorem \cite{OviFur}, the ratio $\D\frac{v_q}{w_q}$ also has the same limit $L$ when $q$ goes to $+\infty$:
\begin{equation}
\D \lim_{q\rightarrow +\infty}\D\frac{v_q}{w_q} = \lim_{q\rightarrow +\infty}\D\frac{\ln\,|u|_{k+q+1,\Omega}}{q} = L,
\end{equation}
and, $|u|_{k+1,\Omega}$ being a constant with respect to $q$,
\begin{equation}\label{exp(-L)}
\D \lim_{q\rightarrow +\infty}\D\left(\frac{|u|_{k+1,\Omega}}{|u|_{k+q+1,\Omega}}\right)^{\frac{1}{q}} = \lim_{q\rightarrow +\infty}\D\left(\frac{1}{|u|_{k+q+1,\Omega}}\right)^{\frac{1}{q}} =e^{-L}=\frac{1}{l}.
\end{equation}
%
As a consequence, from (\ref{Truc}) and (\ref{exp(-L)}) we conclude that:
\begin{equation}
h^*_q \underset{q \to +\infty}{\sim} \frac{1}{e\,l}\,q,
\end{equation}
and
\begin{equation}\label{limite_h*q}
\D \lim_{q\rightarrow +\infty} h^*_q = +\infty.
\end{equation}
$\blacktriangleright$ $\frac{}{}$
Let us now examine the convergence pointwise of the sequence of functions $\D\left(\mathcal{P}_{q}(h)\right)_{q \in \N^\star}$ defined in (\ref{Nonlinear_Prob_Truncate}). \sa
To this end let us, for example, consider a fixed value $h_0$ such that $0 < h_0< h^*_q$. Then, due to (\ref{limite_h*q}), we have:
\begin{equation}
\D\forall\, 0 < h_0 < h^*_q: \lim_{q\rightarrow +\infty}\left(\frac{h_0}{h^*_q}\right)^q = \lim_{q\rightarrow +\infty}e^{\D q\ln\left(\!\frac{h_0}{h^*_q}\right)} = 0^+,
\end{equation}
and similarly for the second part of (\ref{Nonlinear_Prob_Truncate}) corresponding to the case $h_0 > h^*_q$.
\sa
Moreover, when $h_0 = h^*_q, \mathcal{P}_{q}(h^*_q) = \frac{1}{2}, \forall\, q \in\N^*$.
\sa
This enables us to define the pointwise limit function $\mathcal{P}_{0}(h)$ of $\mathcal{P}_{q}(h)$  when $q$ goes to $+\infty$ as:
\begin{equation}\label{Pbar}
\mathcal{P}_{0}(h)= \left |
\begin{array}{ll}
1 & 0 \leq h < +\infty \,,\\
\D\frac{1}{2} & h = +\infty .
\end{array}
\right.
\end{equation}
\sa
Remark that, as $h$ goes to infinity in (\ref{Pbar}), the limit function $\mathcal{P}_{0}(h)$ has a discontinuity.  This comes from the interchange of limiting in $q$ and in $h$ is illicit, namely:
$$
\frac{1}{2}= \lim_{q \rightarrow +\infty} \left[ \lim_{h \rightarrow +\infty} \mathcal{P}_{q}(h)\right] \neq
\lim_{\underset{h \neq h^*_q}{h \rightarrow h^*_q} } \left[  \lim_{q \rightarrow +\infty} \mathcal{P}_{q}(h)\right]   =
\lim_{\underset{h \neq h^*_q}{h \rightarrow h^*_q} } \mathbbm{1}_{]0,+\infty[}(h) =1\,.
$$
\end{prooff}
\vspace{-0.4cm}
\begin{remark}
$\frac{}{}$\\\vspace{-0.4cm}
\begin{itemize}
\item In theorem \ref{Conv_Simple} we assume that the exact solution $u$ to the second order elliptic variational problem \textbf{(VP)} belongs to $H^{r}(\Omega), \forall\, r \in \N$. This is typically the case when the linear form $l(.)$ in (\ref{VP}) is defined by a sufficiently regular function denoted $f$. For example, if $f\in\,H^{r}(\Omega)$ then, for a second order elliptic operator, $u$ belongs to $H^{r+2}(\Omega)$ (see for example \cite{Brezis}).
\item Even if condition (\ref{Cond_ Ratio_Semi_Norm}) of theorem \ref{Conv_Simple} seems restrictive for applications, it is not necessary the case. Take for example the following standard problem:
\begin{equation}\label{edpnum}
\textbf{(VP}\textbf{)} \hspace{0.2cm} \left\{
\begin{array}{l}
\mbox{Find } u \in V \mbox{ solution to:} \\[0.1cm]
-\Delta u = f \mbox{ in } \Omega\,, \\[0.1cm]
u =g \mbox{ on } \partial \Omega\,,
\end{array}
\right.
\end{equation}
where $\Omega$ is the open unit square $]0,1[\times]0,1[$ in $\mathbb{R}^2$, and $f(x,y)=2\pi^2 \sin(\pi x) \cos(\pi y)$. \sa
We readily get that $u(x,y)=\sin(\pi x) \cos(\pi y)$ is the exact solution of (\ref{edpnum}) in  $V\equiv H^{r}(\Omega), (\forall\, r\in\N)$, provided that the Dirichlet boundary condition $g$ is defined by:
\begin{equation}\label{edprhs}
\left\{
\begin{array}{ll}
g(x,0) = \sin(\pi x)   & g(0,y) = 0, \\
g(x,1) = -\sin(\pi x),  & g(1,y) = 0.
\end{array}
\right.
\end{equation}
Then, we obtain that the semi-norm $|u|_{k,\Omega}$ is equal to
\begin{equation}\label{seminorm_exact_sol}
\forall k \geq 0: |u|_{k,\Omega} = (\sqrt{2})^{k-2}\pi^k.
\end{equation}
Finally, on can check that condition (\ref{Cond_ Ratio_Semi_Norm}) is satisfied in that case, as we have
\begin{equation}
\D \lim_{q\rightarrow +\infty}\frac{|u|_{k+q+2,\Omega}}{|u|_{k+q+1,\Omega}} = \frac{1}{\pi\sqrt{2}}.
\end{equation}
\end{itemize}
\end{remark}
\vspace{-0.4cm}
\begin{remark}
$\frac{}{}$\\[-0.5cm]
\begin{itemize}
\item Theorem \ref{Conv_Simple} corresponds to an expected behavior. Indeed, when $q=m-k$ tends to infinity, it claims that the event "$P_m$ \emph{is more accurate that} $P_k$" is an almost sure event for all positive values of $h$. In other words, the higher the distance between $m$ and $k$, the higher the size of the interval $[0,h^*_q]$ where the event "$P_m$ \emph{is more accurate that} $P_k$" is an almost sure event.\vspace{0.2cm}
\item One can notice that this asymptotic feature is also very intuitive in terms of probability. Indeed, as $q=m-k$ goes to infinity,  for $h < h^*_q$, the probability such that $X^{(k+q)}(h)\leq X^{(k)}(h)$ goes to 1, since the domain of existence of $X^{(k+q)}$ goes to 0 whereas the one of $X^{(k)}$ stays fixed and finite. On the contrary, when  $h > h^*_q$ the complementary situation has to be taken into account. Namely, the domain of existence of $X^{(k+q)}$ goes to infinity in comparison with those of $X^{(k)}$ which stays again fixed and finite, (see Figure \ref{Xk+q_vs_Xk}).
\end{itemize}
\end{remark}
\section{Discussion and conclusion}\label{Conclusion}
\noindent In this paper we applied to high order finite elements the novel probabilistic approach we developed in \cite{CMAM_2019} to evaluate the relative accuracy between two Lagrange finite elements $P_k$ and $P_m, (k<m)$. This new way to evaluate the relative accuracy is based on a geometrical interpretation of the error estimate and by considering the approximation errors as random variables. Therefore, we derived two probabilistic laws, the "two steps" one and  the "sigmoid" one, which describe new features of the relative accuracy between finite elements.\sa
%
%
The perspectives of this new approach are not restricted to finite element methods but can be extended to other approximation methods: given a class of numerical schemes and their corresponding error estimates, one is able to order them, not only in terms of asymptotic rate of convergence, but also by evaluating the most probably accurate. \sa
For the finite elements we considered, we can state the following properties as consequences of Theorem \ref{The_nonlinear_law}:
\begin{itemize}
%
\item For the very small values of $h$, the probability such that "$P_m$ \emph{is more accurate than} $P_k$" probability goes to 1. It corresponds to the classical interpretation of Bramble-Hilbert lemma.
\item Depending on the position of $h$ with respect to the critical value $h^*$ defined by (\ref{h*}), $P_k$ or $P_m$ finite elements are more likely accurate.
\item When $h$ is smaller than $h^*$, $P_m$ finite elements are not only asymptotically better than $P_k$ finite elements as $h$ becomes small, but they are \emph{almost surely} more accurate for all of these values of $h$, with a probability between 0.5 and 1.
\item When $h$ is greater than $h^*$, $P_k$ finite elements are \emph{almost surely} more accurate than $P_m$ finite elements, even though $k<m$, with a probability between 0.5 to 1.
\end{itemize}
This last property upsets the widespread idea regarding the relative accuracy between $P_k$ and $P_m, (k<~\!\!m)$, finite elements. It clearly indicates that there exist cases where $P_m$ finite elements \emph{surely} must be overqualified and a significant reduction of implementation time and execution cost can be obtained without a loss of accuracy. We already observed such a phenomenon by using data mining techniques (see \cite{AsCh11}, \cite{AsCh13}, \cite{AsCh16} and \cite{AsCh17}). \sa
However, when the difference between $k$ and $m$ becomes large, high order finite elements are concerned and on can raise the question if $P_k$ finite element is still \emph{almost surely} more accurate than $P_m$ when $h\geq h^*$, or at least, on which interval of $h$ this would still be true. \sa
It is the purpose of the results presented in this paper. To achieve these objectives, we first remark that we need to get asymptotic information regarding the critical value of $h^*_q$ defined in (\ref{h*q_2}) to evaluate the limit of the corresponding sequence of probabilities $\D\left(\frac{}{}\!\!\mathcal{P}_{q}(h)\right)_{q\in\N^\star}$ defined in (\ref{Nonlinear_Prob_Truncate}). \sa
As a consequence, we needed to also get an estimate of the asymptotic behavior of the constant $\mathscr{C}_k$ which is involved in the error estimate (\ref{estimation_error}) derived from Bramble Hilbert Lemma.\sa
So, based on the reference work of Arcangeli and Gout \cite{Arcangeli_Gout}, we identified the necessity to get several estimates which concerned the canonical basis of polynomials $p_i$ introduced in (\ref{shape_function}), related to a given finite element $P_k$ . This is the purpose of Lemma \ref{Estimation_pi}. \sa
Then, we got an explicit $k$-dependency of the constant $\mathscr{C}_k$ in Theorem \ref{C_k_h*q} which enabled us to determine the asymptotic behavior of $h^*_q$, and accordingly, for the sequence of probabilities $\D\left(\frac{}{}\!\!\mathcal{P}_{q}(h)\right)_{q\in\N^\star}$ defined in (\ref{Nonlinear_Prob_Truncate}). This is the purpose of Theorem \ref{Conv_Simple}. \sa
This theorem claims that when one considers two finite elements $P_k$ and $P_m$ for a fixed value of $k$, when $m$ goes to infinity, it does not exist anymore any interval or any value of the mesh size $h$ such that the finite element $P_k$ could be \emph{almost surely} more accurate than $P_m$. \sa
In other words, when one implements high order finite elements $P_m$, the choice of the mesh size $h$ has not to be done to guarantee the better accuracy of the concerned finite element $P_m$ in comparison with another cheaper one $P_k, (k<<m)$.\sa
\textbf{\underline{Acknowledgements}:}
The authors want to warmly dedicate this research to pay homage to the memory of Professors Andr\'e Avez and G\'erard Tronel who largely promote the passion of research and teaching in mathematics.
\end{document}